\documentclass[11pt]{amsart}
\usepackage{amsmath,amsfonts,amssymb}
\usepackage[a4paper]{geometry}
\usepackage{hyperref}
\numberwithin{equation}{section}
\def\d{\mathrm d}
\def\R{\mathbb R}
\def\C{\mathbb C}
\def\D{\mathrm D}
\DeclareMathOperator{\diag}{diag}
\newtheorem{thm}{Theorem}[section]
\newtheorem{lem}[thm]{Lemma}
\newtheorem{cor}[thm]{Corollary}
\newtheorem{prop}[thm]{Proposition}
\theoremstyle{remark}
\newtheorem{expl}{Example}[section]
\newtheorem{rem}{Remark}[section]
\begin{document}
%
\title[Wave equations with variable speed]
{Generalised energy conservation law for wave equations with variable propagation speed}
\author{Fumihiko Hirosawa}
\thanks{Supported by KAKENHI (19740072) from the Ministry of Education, Culture, Sports, Science and Technology (MEXT), Japan.}
\address{Department of Mathematics, Yamaguchi University, 753-8512, Japan}
\email{hirosawa@yamaguchi-u.ac.jp}
\author{Jens Wirth}
\thanks{Supported by EPSRC grant EP/E062873/1 from the Engineering and Physics Research Council, UK.}
\address{Department of Mathematics, Imperial College, London SW7 2AZ, UK}
\email{j.wirth@imperial.ac.uk}
\maketitle
\begin{abstract}
We investigate the long time behaviour of the $L^2$-energy of solutions to wave equations with variable
speed. The novelty of the approach is the combination of estimates for higher order derivatives of
the coefficient with a stabilisation property.  
\end{abstract}
%
\section{Model problem}
We consider the Cauchy problem
\begin{equation}\label{eq:CP}
u_{tt}-a^2(t)\Delta u=0,\qquad u(0,\cdot)=u_1\in H^1(\R^n),\quad \D_t u(0,\cdot)=u_2\in L^2(\R^n)
\end{equation}
for a wave equation with variable propagation speed. As usual we denote $\D_t=-\mathrm i\partial_t$,
$\Delta=\sum_j \partial_{x_j}^2$ the Laplacian on $\R^n$ and $a^2(t)$ is a sufficiently regular non-negative function subject to conditions specified later on. We are interested in the behaviour of the energy as $t\to\infty$ for coefficients bearing {\em very fast oscillations} (in the classification of Reissig-Yagdjian \cite{Yagdjian:2000a}, \cite{Reissig:2004}), but satisfying a suitable {\em stabilisation condition} in the spirit of Hirosawa \cite{Hirosawa:2007}, \cite{Hirosawa:2007b}.  For this we assume that the coefficient $a(t)$ can be written as product 
\begin{equation}
   a(t)=\lambda(t)\omega(t)
\end{equation}
of a shape function $\lambda(t)$ (being essentially free of oscillations) and a bounded perturbation
$\omega(t)$ containing a certain amount of oscillations controlled by our main assumptions. 

Our method leads to an extension of the generalised energy conservation law from \cite{Hirosawa:2007} including the shape function $\lambda(t)$. Roughly speaking, this means that the adapted hyperbolic energy of the solution $u(t,x)$ of \eqref{eq:CP}, 
\begin{equation}
\mathbb E_\lambda(t;u) = \frac12 \int_{\R^n} \big(\lambda^2(t) |\nabla u(t,x)|^2 + |u_t(t,x)|^2\big)\; \d x
\end{equation}
satisfies a {\em two-sided} energy inequality of the form
\begin{equation}
   C_1 \le \frac1{\lambda(t)} \mathbb E_\lambda(t;u) \le C_2
\end{equation}
with constants $C_1$ and $C_2$ depending on the data. The upper bound can be given in terms of
the norms of $u_1\in H^1(\R^n)$ and $u_2\in L^2(\R^n)$, it is {\em not} possible to replace $H^1(\R^n)$ by the corresponding homogeneous space $\dot H^1(\R^n)$ (as in the case of \cite{Hirosawa:2007}).

The behaviour of the energy is only of interest as $t\to\infty$ (or in the neighbourhood of zeros of $\lambda(t)$, which is not within the scope of this note). Therefore it is reasonable to restrict considerations to {\em monotonous} $\lambda(t)$ with $\lambda(0)>0$.

Basic assumptions of our approach are that $a(t)\in C^m(\R_+)$, $m\ge2$, together with
\begin{description}
\item[(A1)] $\lambda(t)>0$, $\lambda'(t)>0$ together with the estimates
\begin{equation}
   \lambda'(t) \approx \lambda(t) \left(\frac{\lambda(t)}{\Lambda(t)}\right),\qquad 
  | \lambda''(t)| \lesssim \lambda(t) \left(\frac{\lambda(t)}{\Lambda(t)}\right)^2,
\end{equation}
where $\Lambda(t) =1+ \int_0^t \lambda(s)\d s$ denotes a primitive of $\lambda(t)$;
\item[(A2)] $0< c_1\le \omega(t)\le c_2$;
\item[(A3)] $\omega(t)$ $\lambda$-stabilises towards 1, i.e. we assume that
\begin{equation}\label{eq:1.5}
   \int_{0}^t \lambda(s) |\omega(s)-1| \d s \lesssim \Theta(t) \ll \Lambda(t), \qquad t\to\infty;
\end{equation} 
\item[(A4)] for $k=1,2,\ldots, m$ the symbol type estimates
\begin{equation}
\left|\d_t^k a(t)\right| \lesssim \lambda(t) \Xi^{-k}(t) 
\end{equation}
are valid, where $\lambda(t)\Xi(t)\gtrsim\Theta(t)$ and
\item[(A5)] 
\begin{equation}
   \int_t^\infty \lambda^{1-m}(s) \Xi^{-m}(s)\d s\lesssim \Theta^{1-m}(t). 
\end{equation}
\end{description}
The number $m$ is determined from (A3)--(A5). The conditions are similar to those from \cite{Hirosawa:2007b}, reason for that is the close relation between wave equations with increasing propagation speed and weakly damped ones. Condition (A5) can be understood as defining property of $\Xi(t)$ in terms of
$\lambda(t)$, the stabilisation rate $\Theta(t)$ and the number $m$.  

Stabilisation condition (A3) is only meaningful if $m\ge 2$. Indeed if (A4) and (A5) hold with $m=1$ 
we would require  $a'(t)/a(t) \in L^1(\R_+)$  and two-sided energy estimates follow directly by Gronwall inequality.

In most examples it is useful to replace assumptions (A4) and (A5) by the following two slightly stronger 
conditions, namely one can use a specific function $\Xi(t)$ depending on $\lambda(t)$, the stabilisation rate $\Theta(t)$ and the number $m$ and assume that
\begin{description}
\item[(A4')] for $k=1,2,\ldots, m$ the symbol type estimates
\begin{equation}
\left|\d_t^k a(t)\right| \lesssim \lambda(t) \left(\frac{\lambda(t)}{\Theta(t)} \left(\frac{\Theta(t)}{\Lambda(t)}\right)^{\frac1m}\right)^k
\end{equation}
are valid and
\item[(A5')] for some number $\epsilon>0$ the estimate $\Lambda^\epsilon(t)\lesssim\Theta(t)$ holds true.
\end{description}
The advantage is that these conditions are more easily checked and the benefit of the number $m$ can be seen directly. Condition (A4') is satisfied for all $m$  if
\begin{description}
\item[(A4'')] for any $\epsilon>0$ and all $k$ the symbol type estimates
\begin{equation}
\left|\d_t^k a(t)\right| \lesssim \lambda(t) \left(\frac{\lambda(t)}{\Theta(t)} \left(\frac{\Theta(t)}{\Lambda(t)}\right)^{\epsilon}\right)^k
\end{equation}
hold true.
\end{description}
Later on we will construct examples along the lines of these conditions and also give counter-examples in the sense that there exists a coefficient satisfying the converse to the inequality (A4'') for $\epsilon<0$ arbitrarily close to $0$ such that the mentioned uniform estimates of the energy do not hold.

{\sl Notational remark:} We use the notation $f\lesssim g$ for two positive functions if there exists a constant
$C$ such that $f\le Cg$ for all values of the arguments. Similarly $f\gtrsim g$ if $g\lesssim f$ and
$f\approx g$ if both $f\lesssim g$ and $g\lesssim f$ are true. Further we denote $f\ll g$ if the quotient
is bounded away from 1, i.e. if $f/g\le c<1$ uniformly in all arguments. For matrices $\|\cdot\|$ denotes
the spectral norm, any other matrix-norm will do as well. Additionally we use $|\cdot|$ for the matrix of the absolute values.

\section{Representation of solutions}
We will not solve \eqref{eq:CP} directly, we will reformulate it as a system of first order and consider the
fundamental solution to that system instead. To be more precise, we apply a partial Fourier transform to reduce \eqref{eq:CP} to an parameter-dependent ordinary differential equation, $\hat u_{tt}+a^2(t)|\xi|^2\hat u=0$, and consider as new unknown the vector
\begin{equation}
V(t,\xi) = \big( \lambda(t)|\xi|\hat u,\D_t \hat u\big)^T.
\end{equation}
We include $\lambda(t)$ to resemble the energy $\mathbb E_\lambda(u;t)=\|V(t,\xi)\|_{L^2}^2$. We could 
include $a(t)$ instead, but in view of (A2) this does not change much. The vector-valued function $V(t,\xi)$ satisfies the first order system
\begin{equation}
\D_t V=  \begin{pmatrix} \frac{\D_t \lambda(t)}{\lambda(t)} & \lambda(t) |\xi| \\ \lambda(t)\omega^2(t)|\xi| & \end{pmatrix} V,
\end{equation}
whose coefficient matrix will be denoted as $A(t,\xi)$. Our aim is to construct the corresponding fundamental solution, i.e. the matrix-valued solution to
\begin{equation}\label{eq:2.3}
  \D_t \mathcal E(t,s,\xi) = A(t,\xi) \mathcal E(t,s,\xi),\qquad \mathcal E(s,s,\xi)=\mathrm I\in \C^{2\times2}.
\end{equation}

If we set formally $\omega(t)=1$ we obtain a much simpler system (by assumption (A1)). Due to its importance for our approach, we denote the corresponding coefficient matrix as $A_\lambda(t,\xi)$
and the corresponding fundamental solution as $\mathcal E_\lambda(t,s,\xi)$. It will be considered first and (partly) constructed in Section~\ref{sec21}.

\subsection{What makes $\lambda(t)$ nice?}\label{sec21}
In a first step we consider the problem with a monotone coefficient. We construct $\mathcal E_\lambda(t,s,\xi)$ for $s,t\ge t_\xi^{(1)}$, where the zone boundary 
$t_\xi^{(1)}$ is given implicitly by 
\begin{equation}\label{eq:2.4}
\Lambda(t_\xi^{(1)})|\xi|=N
\end{equation}
for some fixed constant $N$ (chosen to be sufficiently large) and prove the following statement:

\begin{lem}\label{lem:2.1} Assume (A1). Then
the fundamental solution  $\mathcal E_\lambda(t,s,\xi)$ satisfies uniformly in  $s,t\ge t_\xi^{(1)}$
the two-sided estimate
\begin{equation}
\|\mathcal E_\lambda(t,s,\xi)\| \approx \frac{\sqrt{\lambda(t)}}{\sqrt{\lambda(s)}}
\end{equation}
(regardless of the order of $s$ and $t$).
\end{lem} 

The proof of this fact is essentially given by a $C^2$-theory (in the language of \cite{Hirosawa:2007},
\cite{Hirosawa:2007b})
and follows the corresponding result from \cite{Reissig:2000}.

\begin{pf}
We apply two steps of transformations to the Cauchy problem $\D_t V_\lambda = A_\lambda(t,\xi) V_\lambda$. In  a first one we set $V_\lambda^{(0)}=M^{-1}V_\lambda$, where 
\begin{equation}
  M = \begin{pmatrix} 1&-1\\1&1 \end{pmatrix},\qquad M^{-1} = \frac12  \begin{pmatrix} 1&1\\-1&1 \end{pmatrix}
\end{equation}
is a diagonaliser of the $|\xi|$-homogeneous part of $A_\lambda(t,\xi)$. This yields the new system
\begin{equation}
  \D_t V_\lambda^{(0)} = \left(\begin{pmatrix}\lambda(t)|\xi|&\\&-\lambda(t)|\xi|\end{pmatrix}
  + \frac{\D_t\lambda(t)}{2\lambda(t)}\begin{pmatrix} 1&-1\\-1&1 \end{pmatrix}  \right) V_\lambda^{(0)}.
\end{equation}
For convenience we denote the first (diagonal) matrix as $\mathcal D_\lambda(t,\xi)$ and
the second (remainder) as $R_{0,\lambda}(t,\xi)$. In a second step we want to transform the 
remainder, keeping the structure of the main diagonal part. For this we set
\begin{subequations}
\begin{align}
   N_\lambda(t,\xi) &= \mathrm I + \frac{\D_t\lambda(t)}{4\lambda^2(t)|\xi|} \begin{pmatrix}&1\\-1& \end{pmatrix},\\
   F_\lambda(t,\xi) &=  \frac{\D_t\lambda(t)}{2\lambda(t)}\mathrm I,
\end{align}
\end{subequations}
such that the commutator relation
\begin{equation}
    [\mathcal D_\lambda(t,\xi), N_\lambda(t,\xi)] + R_{\lambda,0}(t,\xi)-F_\lambda(t,\xi) = 0
\end{equation}
holds true. This relation implies that
\begin{align}
  B(t,\xi) &= (\D_t - \mathcal D_\lambda(t,\xi)-R_{\lambda,0}(t,\xi)) N_\lambda(t,\xi) -
  N_\lambda(t,\xi) (\D_t-\mathcal D_\lambda(t,\xi)-F_\lambda(t,\xi)) \notag\\
  &=\D_t N_\lambda(t,\xi) -   [\mathcal D_\lambda(t,\xi), N_\lambda(t,\xi)]  - R_{0,\lambda}(t,\xi) N_\lambda(t,\xi)
  + N_\lambda(t,\xi)F_\lambda(t,\xi) \notag \\
  &= \D_t N_\lambda(t,\xi) -  R_{\lambda,0}(t,\xi) (N_\lambda(t,\xi)-\mathrm I) + (N_\lambda(t,\xi)-\mathrm I)F_\lambda(t,\xi)
\end{align}
is bounded by 
\begin{equation} \label{eq:2.11}
   \| B(t,\xi) \| \lesssim\left| \D_t \frac{\D_t \lambda(t)}{\lambda^2(t)|\xi|} \right|+ \left|\frac{(\D_t\lambda(t))^2}{\lambda^3(t)|\xi|}\right|   \lesssim  \frac{\lambda(t)}{\Lambda^2(t)|\xi|}
\end{equation}
as consequence of assumption (A1). Furthermore, $\|N_\lambda(t,\xi)\|\lesssim 1+ \frac1{\Lambda(t)|\xi|}\lesssim 1$ combined with
\begin{equation}
\det N_\lambda(t,\xi) =  1 -  \frac{(\partial_t \lambda(t))^2}{16\lambda^4(t)|\xi|^2} \ge 1- \frac{C}{N}
\end{equation}
implies that for sufficiently large $N$ the matrix $N_\lambda(t,\xi)$ is invertible with uniformly bounded  inverse, $\|N_\lambda^{-1}(t,\xi)\| \lesssim 1$. This fixes the choice of $N$ for now
(until we may make it slightly larger later on).

Setting $V_\lambda^{(1)}=N_\lambda^{-1}(t,\xi)V_\lambda^{(0)}$ we obtain the system
\begin{equation}\label{eq:2.13}
  \D_t V_\lambda^{(1)} = \left( \mathcal D_\lambda(t,\xi) + F_\lambda(t,\xi) + R_{\lambda,1}(t,\xi) \right) V_\lambda^{(1)}
\end{equation}
with remainder $R_{\lambda,1}(t,\xi) = -N_\lambda^{-1}(t,\xi) B(t,\xi)$ satisfying the bound
\eqref{eq:2.11}. This system can be solved in two steps. First consider the diagonal part 
$\D_t-\mathcal D_\lambda(t,\xi) - F_\lambda(t,\xi)$. The corresponding fundamental solution is
\begin{align}
   \widetilde{\mathcal E}_{\lambda,1}(t,s,\xi) &= \exp\left(\int_s^t (\mathcal D_\lambda(\tau,\xi) + F_\lambda(\tau,\xi))\d\tau\right) \notag\\
   &= \frac{\sqrt{\lambda(t)}}{\sqrt{\lambda(s)}} \diag\left(\mathrm e^{\mathrm i (\Lambda(t)-\Lambda(s))|\xi|},
 \mathrm e^{-\mathrm i (\Lambda(t)-\Lambda(s))|\xi|} \right)
\end{align}
with $\mathrm{cond}\, \widetilde{\mathcal E}_{\lambda,1}(t,s,\xi) = \| \widetilde{\mathcal E}_{\lambda,1}(t,s,\xi) \|\|\widetilde{\mathcal E}_{\lambda,1}(s,t,\xi) \| = 1$. Now, we make the {\sl ansatz}
$\mathcal E_{\lambda,1}(t,s,\xi) = \widetilde{\mathcal E}_{\lambda,1}(t,s,\xi) \mathcal Q_{\lambda,1}
(t,s,\xi)$ for the fundamental solution to \eqref{eq:2.13}. A simple calculation yields for the unknown
$\mathcal Q_{\lambda,1}$ the following equation
\begin{equation}
   \D_t \mathcal Q_{\lambda,1}(t,s,\xi) = \widetilde{\mathcal E}_{\lambda,1}(s,t,\xi) R_{\lambda,1}(t,\xi)\widetilde{\mathcal E}_{\lambda,1}(t,s,\xi) \mathcal Q_{\lambda,1}(t,s,\xi),\qquad  \mathcal Q_{\lambda,1}(s,s,\xi).
\end{equation}
The matrix $\mathcal R_{\lambda,1}(t,s,\xi)= \widetilde{\mathcal E}_{\lambda,1}(s,t,\xi) R_{1,\lambda}(t,\xi)\widetilde{\mathcal E}_{\lambda,1}(t,s,\xi)$ satisfies the bound \eqref{eq:2.11},
\begin{equation}
   \| \mathcal R_{\lambda,1}(t,s,\xi) \| \lesssim \frac{\lambda(t)}{\Lambda^2(t)|\xi|},
\end{equation}
such that the representation of $\mathcal Q_\lambda(t,s,\xi)$ by means of a Peano-Baker series
\begin{equation}
  \mathcal Q_{\lambda,1}(t,s,\xi) = \mathrm I + \sum_{k=1}^\infty
  \int_s^t \mathcal R_{\lambda,1}(t_1,s,\xi) 
  \cdots \int_s^{t_{k-1}}  \mathcal R_{\lambda,1}(t_k,s,\xi) \d t_k\cdots \d t_1
\end{equation}
implies the uniform bound
\begin{align}
  \|  \mathcal Q_{\lambda,1}(t,s,\xi)\| &\le \exp \left( \int_s^t \| \mathcal R_{\lambda,1}(\tau,s,\xi)\| \d\tau\right)
  \le \exp\left( C \int_{t_\xi^{(1)}}^\infty \frac{\lambda(\tau)}{\Lambda^2(\tau)|\xi|}\d\tau\right) \notag\\
 & \le \exp\left(\frac C{\Lambda(t_\xi^{(1)})|\xi|}\right) \le \exp\left(\frac CN\right) \lesssim 1.
\end{align}
The representation $\mathcal E_\lambda(t,s,\xi) = M N_\lambda(t,\xi)  \widetilde{\mathcal E}_{\lambda,1}(t,s,\xi) \mathcal Q_{\lambda,1}N_\lambda^{-1}(s,\xi)M^{-1}$ of the fundamental solution together with the
bounds of all factors established above gives the desired norm estimate. This completes the proof.
$\Box$\end{pf}

\begin{rem}\label{rem:2.1}
In fact we have established more than stated in Lemma~\ref{lem:2.1}. We have a precise description of the structure of the fundamental solution $\mathcal E_\lambda(t,s,\xi)$  which allows to track the large time asymptotics of solutions. To be more precise, we have $N_\lambda(t,\xi)\to\mathrm I$ as $t\to\infty$
for fixed $\xi\ne0$ together with  $\mathcal Q_{\lambda,1}(t,s,\xi) \to  \mathcal Q_{\lambda,1}(\infty,s,\xi)$,
where 
\begin{equation}
   \mathcal Q_{\lambda,1}(\infty,s,\xi) =  \mathrm I + \sum_{k=1}^\infty
  \int_s^\infty  \mathcal R_{\lambda,1}(t_1,s,\xi) 
  \cdots \int_s^{t_{k-1}}  \mathcal R_{\lambda,1}(t_k,s,\xi) \d t_k\cdots \d t_1
\end{equation}
and 
\begin{align}
\|\mathcal Q_{\lambda,1}(t,s,\xi)&-  \mathcal Q_{\lambda,1}(\infty,s,\xi)\|
\le \exp\left(\int_t^\infty \|R_{\lambda,1}(\tau,\xi)\|\d\tau\right)-1\notag\\
&\le \exp\left( \frac C{\Lambda(t)|\xi|}\right) - 1\lesssim \frac C{\Lambda(t)|\xi|}\to 0,\qquad t\to\infty.
\end{align}
Both convergences are locally uniform in $\xi\ne0$. Roughly speaking this means that the solutions 
are determined for large time by $M\widetilde{\mathcal E}_{\lambda,1}(t,s,\xi) M^{-1} $, which is just a free wave (where $\lambda\equiv 1$) with a substitution in the time-variable.
\end{rem}

\subsection{Treatment in the pseudo-differential zone.}  We denote
\begin{equation}
   Z_{pd}(N)=\{(t,\xi)\,:\, 0\le t\le t_\xi^{(1)}\}
\end{equation}
as pseudo-differential zone and continue the construction of the fundamental solution inside this set.
For this we consider $\mathcal E(t,t_\xi^{(1)},\xi)$ and represent its entries as solutions of certain Volterra-type integral equations. Thus, we solve the problem {\em backwards}.

\begin{lem}\label{lem:2.2} Assume (A1) and (A2).
  Then uniformly in $Z_{pd}(N)$ the estimate
  \begin{equation}
     |\mathcal E(t,s,\xi)| \lesssim \begin{pmatrix} \frac{\lambda(t)}{\lambda(s)} & 
     \frac{\lambda(t) (s-t)}{\Lambda(s)}\\ \frac{\lambda(t)}{\lambda(s)} & 1 \end{pmatrix},\qquad t\le s,
  \end{equation}
  holds true.
\end{lem}

Note, that Liouville theorem applied to the original system immediately gives the representation
of the determinant
\begin{equation}
\det\mathcal E(t,s,\xi)=\frac{\lambda(t)}{\lambda(s)},
\end{equation}
which means that we can conclude estimates for the inverse matrix by Cramer's rule.

\begin{cor}\label{cor:2.3}
Uniformly in $Z_{pd}(N)$ the fundamental solution $\mathcal E(t,s,\xi)$ satisfies
\begin{equation}
     |\mathcal E(t,s,\xi)| \lesssim \frac{\lambda(t)}{\lambda(s)} \begin{pmatrix} 1 & 
     \frac{\lambda(s) (t-s)}{\Lambda(t)}\\ \frac{\lambda(s)}{\lambda(t)} & \frac{\lambda(s)}{\lambda(t)} \end{pmatrix}
     \lesssim \begin{pmatrix}\frac{\lambda(t)}{\lambda(s)} & \frac{\lambda(t) (t-s)}{\Lambda(t)} \\ 1 & 1  \end{pmatrix},\qquad s\le t.
  \end{equation}
\end{cor}

\begin{pf}(of Lemma~\ref{lem:2.2})
We consider the columns of $\mathcal E(t,s,\xi)$ seperately and rewrite the differential equation \eqref{eq:2.3} as system of integral equations. This gives for the entries $v(t,\xi)$ and $w(t,\xi)$ of one column
\begin{subequations}
\begin{align}
v(t,\xi) &= \frac{\lambda(t)}{\lambda(s)} v(s,\xi) + \mathrm i |\xi| \lambda(t) \int_s^t w(\tau,\xi)\d\tau, \\
w(t,\xi) &= w(s,\xi) +\mathrm i |\xi|\lambda(t) \int_s^t \omega^2(\tau) v(\tau,\xi) \d\tau
\end{align}
\end{subequations}
with appropriate data $v(s,\xi)$ and $w(s,\xi)$. 

{\sl First column.} We set $v(s,\xi)=1$ and $w(s,\xi)=0$ and restrict to the range $0\le t\le s$. Plugging the second integral equation into the first yields
\begin{align}
v(t,\xi) &= \frac{\lambda(t)}{\lambda(s)} -|\xi|^2 \lambda(t) \int_t^{s} \lambda(\tau) 
\int_\tau^{s} \omega^2(\theta) v(\theta,\xi)\d\theta\d\tau\notag\\
&= \frac{\lambda(t)}{\lambda(s)} -|\xi|^2 \lambda(t)  \int_t^{s} \left( \int_t^\theta \lambda(\tau)\d\tau\right) \omega^2(\theta)v(\theta,\xi)\d\theta.
\end{align} 
The best we can expect is an estimate of the form $\lambda(s)v(t,\xi)/\lambda(t)\in L^\infty(Z'_{pd}(N))$,
where $Z'_{pd}(N)=\{(t,s,\xi)\,:\,0\le t\le s\le t_\xi^{(1)}\}$. Rewriting the integral equation gives 
\begin{equation}
  \frac{\lambda(s)v(t,\xi)}{\lambda(t)} = 1 + \int_t^{s} k_1(t,\theta,\xi) \frac{\lambda(s)v(\theta,\xi)}{\lambda(\theta)} \d\theta
\end{equation}
with kernel
\begin{equation}
  k_1(t,\theta,\xi) =- |\xi|^2 \omega^2(\theta) \lambda(\theta) \int_t^\theta\lambda(\tau) \d\tau, \qquad \theta\in [t,s].
\end{equation}
Now the kernel estimate
\begin{align}
  \sup_{(t,\xi)\in Z_{pd}} \int_0^{s} \sup_{0\le\tilde t\le \theta} |k_1(\tilde t,\theta,\xi)| \d\theta 
  &\lesssim |\xi|^2 \int_0^{t_\xi^{(1)}} \lambda(\theta) \int_0^\theta \lambda(\tau)\d\tau \d\theta\notag\\
  &= |\xi|^2 \int_0^{t_\xi^{(1)}} \Lambda(\theta)\lambda(\theta) \d\theta =\frac12 |\xi|^2\Lambda^2(t_\xi^{(1)})\lesssim 1
\end{align}
uniform in $Z'_{pd}(N)$ implies that the Neumann series 
\begin{equation}
 \frac{\lambda(s)v(t,\xi)}{\lambda(t)} = 1 + \sum_{j=1}^\infty \int_t^{s} k_1(t,t_1,\xi) \cdots \int_{t_{k-1}}^{s} k_1(t,t_k,\xi)\d t_k\cdots \d t_1
\end{equation}
converges in $L^\infty(Z'_{pd}(N))$ (for arbitrary $N$). Therefore, as claimed,
\begin{equation}
|v(t,\xi)| \lesssim \frac{\lambda(t)}{\lambda(s)},
\end{equation}
and the second integral equation implies the corresponding bound for $w(t,\xi)$,
\begin{equation}
|w(t,\xi)| \lesssim |\xi| \lambda(t) \int_t^{s} \frac{\lambda(\tau)}{\lambda(s)}\d\tau 
\le |\xi| \frac{\lambda(t)}{\lambda(s)} \Lambda(t_\xi^{(1)}) \lesssim \frac{\lambda(t)}{\lambda(s)}.
\end{equation}

{\sl Second column.} For the second column we have $v(s,\xi)=0$ and $w(s,\xi)=1$. 
Plugging again the second integral equation into the first one implies
\begin{equation}
  v(t,\xi) = -\mathrm i |\xi| \lambda(t) (s-t) - |\xi|^2 \lambda(t)  \int_t^{s} \left( \int_t^\theta \lambda(\tau)\d\tau\right) \omega^2(\theta)v(\theta,\xi)\d\theta.
\end{equation}
Therefore, we expect $v(t,\xi) / ( |\xi| \lambda(t) (s-t)) \in L^\infty(Z'_{pd}(N))$. Rewriting the
integral equation yields
\begin{equation}
\frac{\mathrm i v(t,\xi)}{ |\xi| \lambda(t) (s-t) } = 1 + \int_t^{s} k_2(t,\theta,\tau) \frac{\mathrm i v(\theta)}{|\xi|\lambda(\theta)(s-\theta)}\d\theta
\end{equation}
with new kernel
\begin{equation}
  k_2(t,\theta,\xi) = - |\xi|^2 \lambda (\theta)\omega^2(\theta) \frac{s-\theta}{s-t} \int_t^\theta \lambda(\tau)\d\tau.
\end{equation}
Note that $|k_2(t,\theta,\xi)|\le |k_1(t,\theta,\xi)|$ (from $t\le\theta\le s$) 
such that the kernel estimate
\begin{equation}
 \sup_{(t,\xi)\in Z_{pd}} \int_0^{s} \sup_{0\le\tilde t\le \theta} |k_2(\tilde t,\theta,\xi)| \d\theta \lesssim 1
\end{equation}
holds true, which in turn implies convergence of the corresponding Neumann series. Therefore, as
claimed,
\begin{equation}
  |v(t,\xi)| \lesssim |\xi| \lambda(t) (s-t) 
\end{equation}
and the second integral equation implies
\begin{equation}
  |w(t,\xi)|  \lesssim 1 + |\xi|^2 \lambda(t) \int_t^{s} \lambda(\theta)(s-\theta)\d\theta 
  \le 1+ |\xi|^2 \lambda(t) \int_t^{t_\xi^{(1)}} (\Lambda(\tau)-1) \d\tau.
\end{equation} 
This is uniformly bounded due to assumption (A1). Indeed, the second term vanishes for
$t=t_\xi^{(1)}$ and its derivative  $$\lambda'(t) \int_t^{t_\xi^{(1)}} (\Lambda(\tau)-1)\d\tau - \lambda(t) (\Lambda(t)-1)$$ 
changes sign. At critical points we get the upper bound $|\xi|^2\lambda^2(t)\Lambda(t) / \lambda'(t)\lesssim |\xi|^2\Lambda^2(t)\lesssim1$ due to the lower bound on $\lambda'(t)$ by (A1).
$\Box$\end{pf}

\subsection{Consideration in the hyperbolic zone.} We define implicitly $t_\xi^{(2)}$ by
\begin{equation}
   \Theta(t_\xi^{(2)})|\xi|=N
\end{equation}
and denote
\begin{equation}
  Z_{hyp}(N) = \{ (t,\xi)\,:\, t\ge t_\xi^{(2)}\}.
\end{equation}
By (A3) we know that $t_\xi^{(2)}>t_\xi^{(1)}$ and $Z_{hyp}(N)$ lies on top of $Z_{pd}(N)$ with a gap
in between. The consideration in the hyperbolic zone follows essentially \cite{Hirosawa:2007} or \cite{Hirosawa:2007b}. Our aim is to obtain the statement of Lemma~\ref{lem:2.1}, but now for the true  $\mathcal E(t,s,\xi)$ and in the smaller zone.

\begin{lem}\label{lem:2.5} Assume (A1), (A2), (A4) and (A5). Then the fundamental solution satisfies 
  \begin{equation}
    \|\mathcal E(t,s,\xi)\| \approx \frac{\sqrt{\lambda(t)}}{\sqrt{\lambda(s)}}
  \end{equation}
  uniformly in $Z_{hyp}(N)$.
\end{lem}

Basically, we follow the proof of Lemma~\ref{lem:2.1}. The main difference is that the remainder terms satisfy worse estimates (due to the presence of $\omega(t)$ in the coefficient matrix), so we do not stop 
after the second step. We apply $m$ steps instead. Before giving the proof we will give this diagonalisation procedure in detail.

We define the following symbol classes within $Z_{hyp}(N)$. We say that $a(t,\xi)$ belongs to
$\mathcal S_N^{\ell}\{m_1,m_2,m_3\}$ if  the symbol estimate
\begin{equation}
   |\D_t^k a(t,\xi)| \le C_k |\xi|^{m_1} \lambda(t)^{m_2} \Xi(t)^{-m_3-k}
\end{equation}
holds true for all $k=0,1,\ldots,\ell$ and all $(t,\xi)\in Z_{hyp}(N)$.
These symbol classes satisfy natural calculus rules. The most important ones for us are collected
in the following proposition.

\begin{prop}
\begin{enumerate}
\item $\mathcal S_N^{\ell}\{m_1,m_2,m_3\}$ is a vector space;
\item $\mathcal S_N^{\ell}\{m_1,m_2,m_3\}\hookrightarrow \mathcal S_{N'}^{\ell'}\{m_1+k,m_2+k,m_3-k\}$
if $N'\ge N$, $\ell'\le\ell$ and $k\ge0$;
\item $\mathcal S_N^{\ell}\{m_1,m_2,m_3\}\cdot \mathcal S_N^{\ell}\{m_1',m_2',m_3'\}\hookrightarrow
\mathcal S_N^{\ell}\{m_1+m_1',m_2+m_2',m_3+m_3'\}$;
\item $\D_t^k\mathcal S_N^{\ell}\{m_1,m_2,m_3\}\hookrightarrow\mathcal S_N^{\ell-k}\{m_1,m_2,m_3+k\}$
for $k\le\ell$;
\item $\mathcal S_N^{0}\{1-m,1-m,m\}\hookrightarrow L^\infty_\xi L^1_t (Z_{hyp}(N))$ with $m$ from assumption (A5).
\end{enumerate}
\end{prop}
Proofs are straightforward. The embedding relation (A2) follows essentially from our requirement $\lambda(t) \Xi(t)\gtrsim \Theta(t)$ in combination with the definition of the zone. 

In order to solve \eqref{eq:2.3} within $Z_{hyp}(N)$ we apply several transformations. In a first step we set $V^{(0)}=M^{-1}(t) V$ with
\begin{equation}
M(t) = \frac1{\omega(t)} \begin{pmatrix} 1& -1 \\ \omega(t) & \omega(t) \end{pmatrix},
\qquad  M^{-1}(t) = \frac12 \begin{pmatrix} \omega(t)& 1 \\ -\omega(t) &1 \end{pmatrix},
\end{equation}
such that 
\begin{equation}
\D_t V^{(0)} = \left( \begin{pmatrix} \lambda(t)\omega(t)|\xi| &\\& -\lambda(t)\omega(t)|\xi| \end{pmatrix}
+ \frac{\D_t (\lambda(t)\omega(t))}{2\lambda(t)\omega(t)} \begin{pmatrix}1&-1\\-1&1 \end{pmatrix}  \right)V^{(0)}
\end{equation}
holds true. Note that the coefficient function $a(t)$ appears in both expressions, such that
the first (diagonal) matrix satisfies $\mathcal D_0(t,\xi)\in\mathcal S_N^{m} \{1,1,0\}$,
while the second (remainder) term is of lower order in our symbol hierarchy $R_0(t,\xi)\in \mathcal S_N^{m-1} \{0,0,1\}$. 

We set $\mathcal D_1(t,\xi)=\mathcal D_0(t,\xi)+\diag R_0(t,\xi)$ and
$R_1(t,\xi)=R_0(t,\xi)-\diag R_0(t,\xi)$. Now we can improve the behaviour of this system within
our symbol classes step by step. 

\begin{lem}\label{lem:2.6}
There exists a zone constant $N$ such that for all $k\le m-1$ we can find matrices
\begin{itemize}
\item $N_k(t,\xi)\in\mathcal S_N^{m-k}\{0,0,0\}$, invertible and 
$N_k^{-1}(t,\xi)\in\mathcal S_N^{m-k}\{0,0,0\}$;
\item $\mathcal D_k(t,\xi)\in\mathcal S_N^{m-k}\{1,1,0\}$ diagonal and\\
$\mathcal D_k(t,\xi)=\diag(\tau_k^+(t,\xi),\tau_k^-(t,\xi))$ with
$|\tau_k^+(t,\xi)-\tau_k^-(t,\xi)|\gtrsim \lambda(t)|\xi|$;
\item $R_k(t,\xi)\in\mathcal S_N^{m-k}\{1-k,1-k,k\}$ antidiagonal  
\end{itemize}
defined on $Z_{hyp}(N)$ such that the operator identity
\begin{equation}
(\D_t-\mathcal D_1(t,\xi)-R_1(t,\xi))N_k(t,\xi) = N_k(t,\xi) (\D_t -\mathcal D_{k+1}(t,\xi)-R_{k+1}(t,\xi))
\end{equation}
holds true.
\end{lem}
\begin{pf}
We construct the matrices $N_k(t,\xi)$ recursively as products
\begin{equation}
  N_k(t,\xi)=\prod_{j=1}^k (\mathrm I+N^{(j)}(t,\xi))
\end{equation}  
of invertible matrices satisfying 
\begin{multline}
(\D_t-\mathcal D_{k}(t,\xi)-R_{k}(t,\xi))(\mathrm I+N^{(k)}(t,\xi)) \\= (\mathrm I+N^{(k)}(t,\xi)) (\D_t -\mathcal D_{k+1}(t,\xi)-R_{k+1}(t,\xi)),\qquad k+1\le m-1.
\end{multline}
This is a straightforward generalisation of the second diagonalisation step in the proof of Lemma~\ref{lem:2.1}. Indeed, the matrices $\mathcal D_1(t,\xi)$ and $R_1(t,\xi)$ satisfy clearly the above
statements. Assume now, the statements about $\mathcal D_k(t,\xi)$ and $R_k(t,\xi)$ 
are true. Then we can construct
\begin{equation}
N^{(k)}(t,\xi) = \begin{pmatrix} & \frac{(R_k(t,\xi))_{12}}{\tau_k^+(t,\xi)-\tau_k^-(t,\xi)} \\
-\frac{(R_k(t,\xi))_{21}}{\tau_k^+(t,\xi)-\tau_k^-(t,\xi)}  &\end{pmatrix} \in\mathcal S_N^{m-k}\{-k,-k,k\},
\end{equation}
such that $\mathrm I+N^{(k)}(t,\xi)$ is invertible for sufficiently large $N$ (following directly from $\|N^{(k)}(t,\xi)\|\lesssim \frac1{|\xi|^k\lambda^k(t)\Xi^k(t)}\lesssim \frac1{|\xi|^k\Theta^k(t)}\le\frac1{N^k}\to 0$ as $N\to\infty$). Furthermore, by construction
\begin{equation}
   [\mathcal D_k(t,\xi),N_k(t,\xi)]+R_k(t,\xi)=0,
\end{equation}
such that
\begin{align}\label{eq:2.50}
    B^{(k)}(t,\xi) &= (\D_t - \mathcal D_k(t,\xi)-R_k(t,\xi)) (\mathrm I+N^{(k)}(t,\xi))
    - (\mathrm I+N^{(k)}(t,\xi)) (\D_t - \mathcal D_k(t,\xi)) \notag\\
    &=\D_t N^{(k)}(t,\xi) - R_k(t,\xi)N^{(k)}(t,\xi)\in\mathcal S_N^{m-k-1}\{-k,-k,k+1\}.
\end{align}
Setting 
\begin{equation}
\mathcal D_{k+1}(t,\xi)=\mathcal D_k(t,\xi) - \diag\big( (\mathrm I+N^{(k)}(t,\xi))^{-1}B^{(k)}(t,\xi)\big)
\end{equation}
and 
\begin{equation}
R_{k+1}(t,\xi) =- (\mathrm I+N^{(k)}(t,\xi))^{-1}B^{(k)}(t,\xi) +  \diag\big( (\mathrm I+N^{(k)}(t,\xi))^{-1}B^{(k)}(t,\xi)\big)
\end{equation}
completes the construction and the symbol estimate of $B^{(k)}$ from \eqref{eq:2.50} finally implies 
$|\tau_{k+1}^+(t,\xi)-\tau_{k+1}^-(t,\xi)|\le |\tau_k^+(t,\xi)-\tau_k^-(t,\xi)|+\lambda(t)|\xi| \frac CN$. 
If we choose $N$ large enough the statement is proven.
$\Box$\end{pf}

\begin{lem}\label{lem:2.7} The diagonal entries satisfy 
\begin{equation}
\mathrm{Im}\,\tau_k^+(t,\xi)=\mathrm{Im}\,\tau_k^-(t,\xi)=-\frac{\lambda'(t)}{2\lambda(t)}-\frac{\omega'(t)}{2\omega(t)}-
\sum_{j=1}^{k-1}  \frac{\partial_t d_j(t,\xi)}{2(d_j(t,\xi)-1)}
\end{equation}
with $d_j(t,\xi)=-\det N^{(j)}(t,\xi)$ being real and $|d_j(t,\xi)|\le c<1$ uniform on $Z_{hyp}(N)$.
\end{lem}
\begin{pf}
The proof goes by induction over $k$. We will show that the above statement and the following hypothesis

\begin{description}
\item[(H$_k$)]
 $R_k(t,\xi)$ has the form $R_k = \mathrm{i}\big(\begin{smallmatrix} & \overline\beta_k \\ \beta_k &\end{smallmatrix}\big)$ with complex-valued $\beta_k(t,\xi)$
\end{description}

are valid. For $k=1$ the assertion (H$_1$) is clearly true with real-valued $\beta_1(t,\xi)=\frac{a'(t)}{2a(t)}$ and $\tau_1^\pm =\pm a(t)|\xi| -\mathrm i \frac{a'(t)}{2a(t)}$ clearly satisfies the statement of Lemma~\ref{lem:2.7}.

We will show that (H$_k$) implies (H$_{k+1}$). The construction implies $N^{(k)} = \frac{\mathrm{i}}{\delta_k} \big(\begin{smallmatrix} & -\overline\beta_k\\\beta_k \end{smallmatrix}\big)$ with $\delta_k(t,\xi)=\tau_k^+(t,\xi)-\tau_k^-(t,\xi)$ being real and 
$|d_k(t,\xi)|=|\det N^{(k)}| = |\beta_k|^2/|\delta_k|^2 \le c < 1$ (for our choice of the zone constant $N$).  Following \cite{Hirosawa:2007} we obtain 
\begin{multline}
 (\mathrm I+N^{(k)})^{-1} (\mathcal D_k+R_k) (\mathrm I+N^{(k)}) 
\\= \frac1{1-d_k} \big( \diag\big(\tau_k^+ - d_k\tau_k^+ - \delta_kd_k, \tau_k^--d_k\tau_k^-+\delta_kd_k\big) +d_kR_k\big) 
\end{multline}
and
\begin{equation}
 (\mathrm I+N^{(k)})^{-1}(\mathrm{D}_t N^{(k)}) = \frac1{1-d_k} \left(\begin{pmatrix} \mathrm{i}\frac{\overline\beta_k}{\delta_k} \partial_t \frac{\beta_k}{\delta_k} &   \\ &
\mathrm{i} \frac{\beta_k}{\delta_k} \partial_t \frac{\overline\beta_k}{\delta_k} 
\end{pmatrix} + \begin{pmatrix}& -\partial_t \frac{\overline\beta_k}{\delta_k}\\\partial_t \frac{\beta_k}{\delta_k} \end{pmatrix}\right)
\end{equation}
such that $\mathrm{Re}\, \frac{\beta_k}{\delta_k} \partial_t \frac{\overline\beta_k}{\delta_k} =\frac{ \partial_td_k}2 = \mathrm{Re}\, 
\frac{\overline\beta_k}{\delta_k} \partial_t \frac{\beta_k}{\delta_k} $ implies
\begin{equation}
\tau_{k+1}^\pm = \tau_k^\pm \mp \frac1{1-d_k} \left(d_k\delta_k +\mathrm{Im}\, \left( \frac{\beta_k}{\delta_k} \partial_t \frac{\overline\beta_k}{\delta_k} 
\right)\right) - \mathrm{i}\frac{\partial_t d_k}{2(d_k-1)}.
\end{equation}
Hence $\delta_{k+1}$ is real again and $R_{k+1}$ satisfies (H$_{k+1}$). Furthermore, the statement
of Lemma~\ref{lem:2.7} follows for $k+1$. 
$\Box$\end{pf}

\begin{pf}(of Lemma~\ref{lem:2.5}) It is sufficient to solve the simpler system 
\begin{equation}
\D_t\mathcal E_m(t,s,\xi)=(\mathcal D_m(t,\xi)+R_m(t,\xi))\mathcal E_m(t,s,\xi),\qquad \mathcal E_m(s,s,\xi)=\mathrm I.
\end{equation}
Lemma~\ref{lem:2.7} implies that the fundamental solution of the diagonal part,
\begin{equation}
\widetilde{\mathcal E}_m(t,s,\xi) = \exp\left(\mathrm i\int_s^t \mathcal D_m(\theta,\xi)\d\theta\right)
=\diag\left(\mathrm e^{\mathrm i\int_s^t \tau_m^+(\theta,\xi)\d\theta},\mathrm e^{\mathrm i\int_s^t \tau_m^-(\theta,\xi)\d\theta}\right),
\end{equation}
has condition number $\mathrm{cond}\,\widetilde{\mathcal E}_m(t,s,\xi)=1$. Therefore, we can make the {\sl ansatz} $\mathcal E_m(t,s,\xi) = \widetilde{\mathcal E}_m(t,s,\xi)\mathcal Q_m(t,s,\xi)$ and get for $\mathcal Q_m(t,s,\xi)$ the system
\begin{equation}
  \D_t \mathcal Q_m(t,s,\xi)=\mathcal R_m(t,s,\xi)\mathcal Q_m(t,s,\xi),\qquad \mathcal Q_m(s,s,\xi)=\mathrm I
\end{equation}
with coefficient matrix $\mathcal R_m(t,s,\xi)=\widetilde{\mathcal E}_m(s,t,\xi)R_m(t,\xi)\widetilde{\mathcal E}_m(t,s,\xi)$ subject to the same bounds like $R_m(t,\xi)$, 
\begin{equation}
\|\mathcal R_m(t,s,\xi)\|=\|R_m(t,\xi)\|\lesssim \frac 1{|\xi|^{1-m}\lambda^{1-m}(t)\Xi^m(t)}.
\end{equation}
Therefore, $\mathcal Q_m(t,s,\xi)$ can be represented as Peano-Baker series and satisfies the
uniform estimate
\begin{align}\label{eq:2.57}
  \| \mathcal Q_m(t,s,\xi)\|&\le \exp\left(\int_s^t \|\mathcal R_m(\theta,s,\xi)\|\d\theta \right)
  \le \exp\left(\int_{t_\xi^{(2)}}^\infty  \frac C{|\xi|^{1-m}\lambda^{1-m}(\theta)\Xi^m(\theta)}\d\theta\right)\notag\\
  &\le \exp\left( \frac C{|\xi|^{1-m}\Theta^{1-m}(t_\xi^{(2)})} \right) \lesssim 1.
\end{align}
Additionally, by Liouville theorem and the invariance of the trace under similarity transformations
we get
\begin{equation}
 \det\mathcal Q_m(t,s,\xi) = \exp\left(\mathrm i\int_s^t \mathrm{tr}\, \mathcal R_m(\theta,s,\xi)\d\theta\right)
 =  \exp\left(\mathrm i\int_s^t \mathrm{tr}\, R_m(\theta,\xi)\d\theta\right) = 1
\end{equation}
and $\|\mathcal Q_m^{-1}(t,s,\xi)\|\lesssim1$.
Thus, representing $\mathcal E(t,s,\xi)$ as
\begin{equation}
\mathcal E(t,s,\xi) = M(t) N_m(t,\xi)  \widetilde{\mathcal E}_m(t,s,\xi)\mathcal Q_m(t,s,\xi) N_m^{-1}(s,\xi)M^{-1}(s)
\end{equation}
gives by the uniform bounds of \eqref{eq:2.57} and Lemma~\ref{lem:2.6} 
\begin{equation}
  \| \mathcal E(t,s,\xi) \| \approx \|\widetilde{\mathcal E}_m(t,s,\xi)\| 
  = \exp\left(-\int_s^t \mathrm{Im}\, \tau_m^\pm(\theta,\xi)\d\theta\right) 
  \approx \frac{\sqrt{\lambda(t)}}{\sqrt{\lambda(s)}}
\end{equation}
and the statement is proven.
$\Box$\end{pf}

\begin{rem}\label{rem:2.2}
Again, we have established much more than just the two-sided estimate of Lemma~\ref{lem:2.5}. We
got a precise description of the structure of the fundamental solution $\mathcal E(t,s,\xi)$ for large
time $t$. Indeed, like in Remark~\ref{rem:2.1} about Lemma~\ref{lem:2.1} we established that the
transformation matrices $N_k(t,\xi)\to\mathrm I$ and the amplitudes
$\mathcal Q_m(t,s,\xi)\to\mathcal Q_m(\infty,s,\xi)$ as $t\to\infty$ locally uniform in $\xi\ne0$.
Therefore, solutions are determined for large time by $M(t)\widetilde{\mathcal E}_m(t,s,\xi)$.
\end{rem}

\subsection{Consideration in the intermediate zone} This zone is defined as 
\begin{equation}
  Z_{int}(N) = \{ (t,\xi)\,:\, t_\xi^{(1)}\le t\le t_\xi^{(2)}\}.
\end{equation}
We want to relate $\mathcal E(t,s,\xi)$ to $\mathcal E_\lambda(t,s,\xi)$ within this zone. For this we employ the stabilisation condition {\em in combination} with Lemma~\ref{lem:2.1}.

\begin{lem}\label{lem:2.8} Assume (A1) -- (A3). Then the fundamental solution satisfies
  \begin{equation}
    \|\mathcal E(t,s,\xi)\| \approx \frac{\sqrt{\lambda(t)}}{\sqrt{\lambda(s)}}
  \end{equation}
  uniformly in $Z_{int}(N)$.
\end{lem}
\begin{pf}
We make the {\em ansatz} $\mathcal E(t,s,\xi)=\mathcal E_\lambda(t,s,\xi)\mathcal Q_{int}(t,s,\xi)$. Then
the matrix $\mathcal Q_{int}(t,s,\xi)$ satisfies the differential equation
\begin{equation}
\D_t \mathcal Q_{int}(t,s,\xi)=\mathcal E_\lambda(s,t,\xi)\big(A(t,\xi)-A_\lambda(t,\xi)\big) \mathcal E_\lambda(t,s,\xi) \mathcal Q_{int}(t,s,\xi)
\end{equation} 
with initial condition $\mathcal Q_{int}(s,s,\xi)=\mathrm I$. The stabilisation condition together with
the uniform bound of the condition number $\mathrm{cond}\,\mathcal E_\lambda(t,s,\xi)\lesssim1$ 
from Lemma~\ref{lem:2.1} implies that the coefficient matrix of this problem satisfies 
\begin{align}
   &\int_{t_\xi^{(1)}}^{t_\xi^{(2)}} \|\mathcal E_\lambda(s,t,\xi)\big(A(t,\xi)-A_\lambda(t,\xi)\big) \mathcal E_\lambda(t,s,\xi)\| \d t\notag \\
   \lesssim & \int_{t_\xi^{(1)}}^{t_\xi^{(2)}} \|A(t,\xi)-A_\lambda(t,\xi)\| \d t\notag \\
   \approx & |\xi|  \int_{t_\xi^{(1)}}^{t_\xi^{(2)}} \lambda(t) |\omega^2(t)-1|\d t \lesssim |\xi|\Theta(t_\xi^{(2)})=N.
\end{align}
Therefore, the representation of $\mathcal Q_{int}(t,s,\xi)$ as Peano-Baker series implies the uniform boundedness of $\mathcal Q_{int}(t,s,\xi)$ over the  intermediate zone. Furthermore, we get $\det\mathcal Q_{int}(t,s,\xi)=1$ from Liouville theorem and conclude that $\mathcal Q_{int}(t,s,\xi)$ is uniformly invertible. This transfers the two-sided estimate from $\mathcal E_\lambda(t,s,\xi)$ to $\mathcal E(t,s,\xi)$
and the statement is proven.
$\Box$\end{pf}

\section{Energy inequalities}
\subsection{Estimates from above} The statements of Lemmata~\ref{lem:2.5} and \ref{lem:2.8} imply that the energy $\mathbb E_\lambda(u;t)$ increases (for large $t$ and $\xi$) like $\lambda(t)$. Our first aim
is to combine this with the estimate from Lemma~\ref{lem:2.2} / Corollary~\ref{cor:2.3}. For this we assume that 
\begin{description}
\item[(A1$_+$)] the coefficient function $\lambda(t)$ satisfies $t\sqrt{\lambda(t)}\lesssim\Lambda(t)$ in addition to (A1),
\end{description}
which is true in all example cases. This might even be a consequence of (A1), however we don't know that for certain.
%


\begin{thm}\label{thm:3.2} Assume (A1$_+$) -- (A5). Then 
all solutions $u(t,x)$ to the Cauchy problem \eqref{eq:CP} satisfy the {\sl a priori} estimate  
\begin{equation}
	 \mathbb E_\lambda(u;t) \le C \lambda(t) \left(\|u_1\|_{H^1}^2+\|u_2\|_{L^2}^2\right)
\end{equation}
with a constant $C$ depending only on the coefficient function $a(t)$.
\end{thm}
\begin{pf}
Corollary~\ref{cor:2.3} implies the estimate
\begin{equation}
   \| \mathcal E(t,0,\xi)\,\diag(|\xi|/\langle\xi\rangle,1)\| \lesssim \max(t|\xi|\lambda(t),1)
   \lesssim \sqrt{\lambda(t)}
\end{equation}
in combination with (A1$_+$) and the definition of the pseudo-differential zone. By Lemmata~\ref{lem:2.5} and \ref{lem:2.8}
\begin{align}
   \| \mathcal E(t,0,\xi)\,\diag(|\xi|/\langle\xi\rangle,1)\|&\lesssim 
   \|\mathcal E(t,t_\xi^{(1)},\xi)\|\,\| \mathcal E(t_\xi^{(1)},0,\xi)\,\diag(|\xi|/\langle\xi\rangle,1)\|\notag\\
    &\lesssim  \sqrt{\lambda(t)}.
\end{align}
follows for all $t\ge t_\xi^{(1)}$ and the proof is complete. 
$\Box$\end{pf}

\begin{rem}
If $\lambda(t)$ is bounded we do not need to change the space for the data. The additional factor $|\xi|$ for small frequencies used in the previous argument to compensate the estimate of Corollary~\ref{cor:2.3} is not necessary in this case and the statement
\begin{equation}
   \mathbb E_\lambda(u;t)\lesssim \mathbb E_\lambda(u;0),\qquad \lambda(t)\le c<\infty,
\end{equation}
from \cite{Hirosawa:2007} follows. However, if $\lambda(t)$ is unbounded, this estimate is in general false. This can be seen by constructing explicit representations in terms of special functions (like done for 
$a(t)=t^\ell$ in \cite{Reissig:1997b} or $a(t)=\mathrm e^t$ in \cite{Galstian:2003}) and evaluating them in the neighbourhood of $\xi=0$. See also \cite{Wirth:2004c} for a similar argument in the dissipative case.
\end{rem}

\subsection{Bounds from below} 
Outside the pseudo-differential zone we already achieved lower bounds. Our strategy is to relate solutions to a quantity which can be controlled everywhere. This idea will be combined with an application of Banach-Steinhaus theorem on a dense subspace of $H^1(\R^n)\times L^2(\R^n)$ {\em excluding} the exceptional frequency $\xi=0$. 

\begin{thm} Assume (A1$_+$) -- (A5). Then for all data $u_1\in H^1(\R^n)$ and $u_2\in L^2(\R^n)$ there exists a constant $C$ such
that 
\begin{equation}
  \mathbb E_\lambda(u;t)\ge C{\lambda(t)} 
\end{equation}
holds true. (The constant $C$ depends in a nontrivial way on the data.)
\end{thm}
\begin{pf}
We proceed in two steps. In a first step we assume that the data $u_1$ and $u_2$ satisfy
the condition $0\not\in\mathrm{supp}\,\hat u_i$ (which directly implies $0\not\in\mathrm{supp}\,\hat u(t,\cdot)$ for all $t\ge0$). We want to compare $(|\xi|\hat u,\D_t\hat u)^T$ to $\widehat{\mathcal E}(t,0,\xi) (\tilde w_1,\tilde w_2)^T$, where
\begin{equation}
   \widehat{\mathcal E}(t,0,\xi) =\begin{cases} \sqrt{\lambda(t_\xi^{(1)})}\, \mathcal E(t,t_\xi^{(1)},\xi),\qquad &t\ge t_\xi^{(1)},\\ \sqrt{\lambda(t)}\,\mathrm I,&0\le t\le t_\xi^{(1)}, \end{cases}
\end{equation}
for suitably chosen $\tilde w_i\in L^2(\R^n)$. By definition and Lemmata~\ref{lem:2.5} and \ref{lem:2.8}
we have $\|  \widehat{\mathcal E}(t,0,\xi) \|\approx \sqrt{\lambda(t)}$ and $\| \widehat{\mathcal E}^{-1}(t,0,\xi) \|\approx 1/\sqrt{\lambda(t)}$ such that the two-sided estimate
\begin{equation}\label{eq:3.8}
   \| \widehat{\mathcal E}(t,0,\xi) (\tilde w_1,\tilde w_2)^T \|_2 \approx \sqrt{\lambda(t)} \| (\tilde w_1,\tilde w_2)\|_2
\end{equation}   
follows. Now we will construct $\tilde w_i$ such that
\begin{equation}\label{eq:3.9}
  \frac1{\sqrt{\lambda(t)}}\| (|\xi|\hat u,\D_t\hat u)^T - \widehat{\mathcal E}(t,0,\xi) (\tilde w_1,\tilde w_2)^T\|_2\to 0,\qquad t\to\infty.
\end{equation}
Since $0\not\in\mathrm{supp}\,\hat u(t,\cdot)$ the difference vanishes identically for sufficiently large $t$ if we define
\begin{align}\label{eq:3.11}
 (\tilde w_1,\tilde w_2)^T &= \lim_{t\to\infty}  \widehat{\mathcal E}^{-1}(t,0,\xi)
\mathcal E(t,0,\xi) \diag(|\xi|/\langle\xi\rangle,1) (\langle\xi\rangle\hat u_1,\hat u_2)^T\notag\\
&=\frac1{\sqrt{\lambda(t_\xi^{(1)})}} \mathcal E(t_\xi^{(1)},0,\xi) \diag(|\xi|/\langle\xi\rangle,1) (\langle\xi\rangle\hat u_1,\hat u_2)^T
\end{align}
and by the argument used in the previous proof the appearing multiplier is uniformly bounded in $\xi$.
Thus for all data with  $0\not\in\mathrm{supp}\,\hat u_i$ we constructed $\tilde w_i\in L^2(\R^n)$.

In a second step we relax the condition on the data. This follows by Banach-Steinhaus theorem
since we are already on a dense subset of $[L^2(\R^n)]^2$ and the left hand side of \eqref{eq:3.9}
is uniformly bounded by Theorem~\ref{thm:3.2}. Thus \eqref{eq:3.9} holds for all solutions if we define
$\tilde w_i$ by \eqref{eq:3.11} in terms of the data.

Finally, from \eqref{eq:3.8} and \eqref{eq:3.9} the desired statement follows.
$\Box$\end{pf}

\section{Examples and counter-examples}
We will collect some examples for shape functions $\lambda(t)$ and perturbations $\omega(t)$ which are admissible in our context. At first we introduce several classes examples depending on the growth order of $\lambda(t)$ and give suitable $\Theta(t)$ and $\Xi(t)$ for assumptions (A1) to (A5). Later on we construct functions $\omega(t)$ subject to corresponding the bounds in all these cases. 

Finally Section~\ref{sec4.3} is devoted to counter-examples, i.e. to show that the symbol-type assumption (A4'') for the coefficient is indeed sharp within certain classes of examples.

\subsection{Classes of examples}

\begin{expl}\label{expl1}
(Polynomial growth) It is possible to choose all functions as polynomials. To be precise, we can set
\begin{subequations}
\begin{align}
  \lambda(t) &= (1+t)^p,\\ 
  \Theta(t)&=(1+t)^{1+q},\\
  \Xi(t)&=(1+t)^r
\end{align}
for suitable choices of $p$, $q$ and $r$. For any $p>0$ assumption (A1) is fulfilled. Furthermore, we need
$0\le q<p$ for (A3) and 
\begin{equation}
\begin{cases}
  1\ge r \ge  r_m=1-p+q+\frac{p-q}m ,\qquad & \text{for (A4') and (A5')},\\
  1\ge r > r_\infty = 1-p+q, &\text{for (A4'') and (A5')}.
\end{cases}
\end{equation}
Increasing $m$ makes $r_m$ smaller and therefore the symbol condition (A4') becomes weaker for fixed derivatives (however, we need more derivatives). In this sense stabilisation allows to weaken symbol estimates. 
\end{subequations}
\end{expl}

\begin{expl}\label{expl2}
(Suprapolynomial growth) It is of interest to look at problems with faster increasing $\lambda(t)$. Therefore, we consider
\begin{subequations}
\begin{align}
  \lambda(t) &= \exp(t^\alpha),\qquad\alpha\in(0,1)\\ 
  \Theta(t)&=t^{-\beta} \exp(t^\alpha),\\
  \Xi(t)&=t^\gamma.
\end{align}
Again we check all the requirements. Assumption (A1) is fulfilled. For (A3) we need $\beta>\alpha-1$ and 
\begin{equation}
\begin{cases}
  1-\alpha\ge \gamma\ge\gamma_m = -\beta+\frac{\beta-\alpha+1}m,\qquad & \text{for (A4') and (A5')},\\
    1-\alpha\ge \gamma>\gamma_\infty= -\beta, &\text{for (A4'') and (A5')}.
\end{cases}
\end{equation}
Again increasing $m$ decreases $\gamma_m$ and the interesting values for $\gamma$ are negative.
\end{subequations}
\end{expl}

\begin{expl}\label{expl3}
(Exponential growth) It is not essential that $\Xi(t)$ is polynomial. We can also consider
\begin{subequations}
\begin{align}
  \lambda(t) &=\mathrm e^t\\ 
  \Theta(t)&=\mathrm e^{at},\\
  \Xi(t)&=\mathrm e^{bt}
\end{align}
under suitable conditions on $a$ and $b$. Assumption (A1) is fulfilled. For (A3) we need 
$a<1$ and 
\begin{equation}
\begin{cases}
  0\ge b\ge b_m = a-1 + \frac{1-a}m,\qquad & \text{for (A4') and (A5')},\\
  0\ge b>b_\infty= a-1, &\text{for (A4'') and (A5')}.
\end{cases}
\end{equation}
\end{subequations}
\end{expl}

\subsection{Construction of admissible $\omega(t)$} Nontrivial examples for perturbations
$\omega(t)$ of the `nice' coefficient $\lambda(t)$ can be constructed in all cases. Our method depends on the choice of three positive sequences, 
\begin{equation}\label{eq:4.4}
t_j\to\infty,\qquad \delta_j\le\Delta t_j=t_{j+1}-t_j\quad\text{ and }\quad\eta_j\le 1
\end{equation}
and  a function $\psi\in C_0^m(\R)$ with 
\begin{equation}\label{eq:4.5}
\mathrm{supp}\,\psi\subseteq[0,1],\quad -1<\psi(t)<1 \quad\text{and}\quad \int_0^1|\psi(t)|\d t =\frac12.
\end{equation}
Using these ingredients we define
\begin{equation}\label{eq:4.6}
\omega(t) = 1 + \sum_{j=1}^\infty  \eta_j \psi\left(\frac{t-t_j}{\delta_j}\right), 
\end{equation}
the sum is converging trivially, since for each $t$ at most one term is present. Furthermore, if $c_1=\min\psi(t)$ and $c_2=\max\psi(t)$  then we get the bound $0<1+c_1 \le \omega(t) \le 1+c_2$.  It remains to look at the stabilisation properties and the symbol estimates. For the first one note that
\begin{equation}
  \int_0^t \lambda(s) |\omega(s)-1|\d s = \sum_{j=1}^k \eta_j \int_{t_{j}}^{t_{j+1}} \lambda(s) \left|\psi\left(\frac{s-t_j}{\delta_j}\right)\right|\d s
  \le  \sum_{j=1}^k \eta_j \delta_j \lambda(t_{j+1})
\end{equation}
for $t\in[t_{k},t_{k+1}]$. Similarly, we get the lower bound 
$\sum_{j=1}^k \eta_j \delta_j \lambda(t_{j})$. Stabilisation property (A3) is ensured, if $\eta_j\delta_j$ are small enough to guarantee 
\begin{equation}\label{eq:4.8}
  \Theta(t_{k+1})\approx\sum_{j=1}^k \eta_j \delta_j \lambda(t_{j+1}) \ll \sum_{j=1}^k \lambda(t_j) \Delta t_j \le \Lambda(t_{k+1}).
\end{equation}
Derivatives of $\omega(t)$ can be estimated by a multiplication with $\delta_j^{-1}$ on $[t_j,t_{j+1}]$,
such that $\Xi(t)$ should satisfy $\Xi(t_j)\lesssim \delta_j$.

\begin{expl}(Polynomial case) 
We consider $\lambda(t)=(1+t)^p$ from Example~\ref{expl1} and give a suitable choice of sequences. We choose $t_j = 2^j$, such that $\Delta t_j=2^{j-1}$ and parameters $p$, $q$ and $r$ from Example~\ref{expl1}. Then $\delta_j$ is determined by $\delta_j \approx \Xi(t_j)$ as $\delta_j= 2^{jr-1}$ and \eqref{eq:4.8} implies our choice for $\eta_j$,
\begin{equation}
   \eta_j=2^{j(1+q-p-r)}.
\end{equation} 
Due to $r\ge r_m=1+p-q+(p-q)/m$ this choice implies $0<\eta_j\le 1$. 
\end{expl}

\begin{expl}(Suprapolynomial case)
We consider $\lambda(t)=\exp(t^\alpha)$ from Example~\ref{expl2}. To simplify the summation in \eqref{eq:4.8} we adjust $t_j$ such that $\lambda(t_j)\approx \mathrm e^j$. This gives $t_j=j^{1/\alpha}$,
$\Delta t_j\ge \frac1\alpha j^{1/\alpha-1}$. We choose $\delta_j=j^{\gamma/\alpha}$ (which is smaller than $\Delta t_j$ due to $\gamma<1-\alpha$) and $\eta_j=j^{-(\beta+\gamma)/\alpha}$, such that the left part of \eqref{eq:4.8} is satisfied.
\end{expl}

\begin{expl}(Exponential case)
We consider $\lambda(t)=\mathrm e^t$ from Example~\ref{expl3}. In this situation we choose
$t_j=j$ and determine the sequences in dependence of the given parameters $a$ and $b$ from Example~\ref{expl3}. This implies $\delta_j=\mathrm e^{bj}$ and $\eta_j=\mathrm e^{j(a-b-1)}$.
By assumption $b<0$ and $a-b-1\le0$ and therefore $\delta_j<1$ and $\eta_j\le 1$.
\end{expl}

\subsection{Counter-examples}\label{sec4.3} Finally we want to apply a modified Floquet approach to show that our considerations are optimal in the sense that for given $\lambda(t)$ from our example classes there exists a coefficient $\omega(t)$ which violates one of the assumptions nearly and in turn leads to the non-existence of uniform bounds. The approach is a generalisation of considerations from \cite{Reissig:1999}, \cite{Tarama:1995} and implicitly also used in \cite{Hirosawa:2007}.

The construction of the coefficient function follows that from the previous section with one alteration, we do not just add one bump $\psi(t)$ in the intervals $[t_j,t_{j+1}]$ but $\nu_j$ many of them.
Thus we are given sequences $t_j$, $\delta_j$ subject to \eqref{eq:4.4} and $\nu_j\in\mathbb N$  together with a real-valued  function $\psi\in C_0^\infty[0,1]$ subject to \eqref{eq:4.5} and 1-periodised as $b(t)=\psi(t \mod 1)$. Then $\omega(t)$ is given by
\begin{equation}\label{eq:4.6}
\omega(t) =\begin{cases} 1, & t\not\in\bigcup_{j=1}^\infty [t_j,t_j+\delta_j],\\
 1 +  b\left(\frac{\nu_j}{\delta_j}(t-t_j)\right), & t\in  [t_j,t_j+\delta_j] .
 \end{cases}
\end{equation}
All parameters are adjusted in a suitable way in dependence of the given $\lambda(t)$. Stabilisation is guaranteed if $\Theta(t_{k+1})\approx \sum_{j=1}^k \delta_j\lambda(t_{j+1})$ is small compared to $\Lambda(t_{k+1})$ and derivatives behave like multiplication with $\nu_k/\delta_k$ on $[t_k,t_k+\delta_k]$, i.e. we have to impose $\Xi(t)\lesssim \delta_k/\nu_k$ for $t\in[t_k,t_k+\delta_k]$.  By adjusting the sequence $\delta_j$ we can influence the stabilisation rate, while adjusting $\nu_j$ allows to change the symbolic estimates.

\subsubsection{A lower estimate for the fundamental solution on $[t_j,t_j+\delta_j]$}  We introduce a new local  time-variable $s$ such that $t(s)=t_j+s\delta_j/\nu_j$, $s\ge0$, and look for the fundamental solution
$\mathcal Y_j(s,s_0,\xi):=\mathcal E(t(s),t(s_0),\xi)$.  This matrix-valued function satisfies
\begin{equation}
   \D_s \mathcal Y_j(s,s_0,\xi) =A_j(s,\xi)\mathcal Y_j(s,s_0,\xi),\qquad \mathcal Y_j(s_0,s_0,\xi)=\mathrm I
\end{equation}
with coefficient matrix
\begin{equation}
A_j(s,\xi)=  \frac{\delta_j}{\nu_j} A(t(s),\xi) =\frac{\delta_j}{\nu_j}  \begin{pmatrix} -\mathrm i\frac{\lambda'(t(s))}{\lambda(t(s))}& \lambda(t(s))|\xi| \\
  \lambda(t(s)) (1+b(s))^2 |\xi| \end{pmatrix},\quad s\in[0,\nu_j].
\end{equation}
Our strategy is to relate this to the $j$-independent periodic problem with coefficient matrix
\begin{equation}
  B(s,\tilde\lambda) = \begin{pmatrix} & \tilde \lambda  \\
  \tilde\lambda\,  (1+b(s))^2  \end{pmatrix},\qquad  s\in\R,
\end{equation}
and parameter $\tilde\lambda= \delta_j \lambda(t_j)|\xi| / \nu_j$, i.e. to consider
\begin{equation}\label{eq:Hill}
   \D_s \mathcal X(s,\tilde\lambda) = B(s,\tilde\lambda) \mathcal X(s,\tilde\lambda), \qquad \mathcal X(0,\tilde\lambda)=\mathrm I.
\end{equation} 
Periodicity of the problem allows to restrict most considerations to the monodromy matrix $\mathcal X(\tilde\lambda)=\mathcal X(1,\tilde\lambda)$. An elementary application of Floquet theory (based on 
$1+b(s)$ strictly positive, $b\in C^2(\R)$, $1$-periodic and real-valued) implies

\begin{lem}[Floquet theorem, cf. \cite{Magnus:1966}]\label{lem:4.1}
There exists a bounded open subinterval $\mathcal I$ of $(0,\infty)$ such that the monodromy matrix $\mathcal X(\tilde\lambda)$ of \eqref{eq:Hill} has for all parameters $\tilde\lambda\in\mathcal I$ a purely imaginary eigenvalue of magnitude larger than 1.
\end{lem}

Thus in order to get the worst possible behaviour of solutions we restrict our considerations to 
\begin{equation}
  \xi\in\Omega_j := \{ \xi\in\R^n\;:\;\tilde\lambda=\frac{\delta_j\lambda(t_j)}{\nu_j} |\xi|\in\mathcal I\}.
\end{equation}
It is evident that $\Omega_j$ is of positive measure, even if we shrink $\mathcal I$ in such a way that we have a uniform lower bound for the magnitude of the eigenvalue. We use Lemma~\ref{lem:4.1} to show that the following statement holds true for $\mathcal Y_j(\nu_j,0,\xi)$ uniform in $j$ and $\xi\in\Omega_j$.

\begin{lem}\label{lem:4.2}
Assume $\delta_j\frac{\lambda(t_j)}{\Lambda(t_j)}\to0$, $\lambda(t_j+\delta_j) \approx \lambda(t_j)$ and
$\Lambda(t_j+\delta_j)\approx \Lambda(t_j)$  uniform in $j$. 
Then there exists $\mu>1$ depending on $b(s)$ and the choice of $\mathcal I$, such  the matrix $\mathcal Y_j(\nu_j,0,\xi)$ has for all $\xi\in\Omega_j$ and sufficiently large $j$ an eigenvalue of modulus greater than $\frac{\mu^{\nu_j}}2$.
\end{lem}
\begin{pf} {\sl Step 1.}
We write $\mathcal Y_j(\nu_j,0,\xi)=\mathcal Y_j(\nu_j,\nu_j-1,\xi)\cdots\mathcal Y_j(2,1,\xi) \mathcal Y_j(1,0,\xi)$ and prove the estimates
\begin{align}
&\|\mathcal Y_j(k+1,k,\xi) - \mathcal Y_j(k,k-1,\xi)\|\lesssim \frac{\delta_j}{\nu_j}\frac{\lambda(t_j)}{\Lambda(t_j)} ,\label{eq:4.16}\\
&\| \mathcal Y_j(k+1,k,\xi) - \mathcal X(\tilde\lambda )\|\lesssim  \delta_j\frac{\lambda(t_j)}{\Lambda(t_j)} ,
\qquad\tilde\lambda= \delta_j\lambda(t_j)|\xi|/{\nu_j},\label{eq:4.17}
\end{align}
for $k=(0,)1,\cdots \nu_j-1$ uniform in $\xi\in\Omega_j$ and $j$. Note for this, that uniform in $j$, $\tau\in[0,1]$ and $k$ in the above stated ranges 
\begin{align}
   &\|A_j(k+\tau,\xi)-A_j(k+\tau-1,\xi)\| \approx \frac{\delta_j}{\nu_j} |\lambda(t(k+\tau))-\lambda(t(k-1+\tau))| |\xi|\notag\\& + \frac{\delta_j}{\nu_j}\left| \frac{\lambda'(t(k+\tau)}{\lambda(t(k+\tau)} -\frac{\lambda'(t(k-1+\tau)}{\lambda(t(k-1+\tau)} \right|  \lesssim \frac{\delta_j^2}{\nu_j^2} |\xi| \lambda'(t(\zeta))
+ \frac{\delta_j^2}{\nu_j^2} \left(\frac{\lambda(t(\zeta))}{\Lambda(t(\zeta))}\right)^2\notag\\&    \lesssim \frac{\delta_j^2}{\nu_j^2} |\xi|\frac{\lambda^2(t_j)}{\Lambda(t_j)}\approx \frac{\delta_j}{\nu_j}\frac{\lambda(t_j)}{\Lambda(t_j)}  \\
&  \|A_j(k+\tau,\xi)-B(\tau, \tilde\lambda) \| \approx  \frac{\delta_j}{\nu_j} |\lambda(t(k+\tau))-\lambda(t(0))| |\xi| +  \frac{\delta_j}{\nu_j}\frac{\lambda'(t(k+\tau))}{\lambda(t(k+\tau))} \notag\\& \lesssim 
  \frac{\delta_j^2}{\nu_j^2} |\xi| (k \lambda'(t(\zeta))+\lambda'(t(k+\tau)))\lesssim \delta_j \frac{\lambda(t_j)}{\Lambda(t_j)}    ,  
\end{align}
hold true (with intermediate values $\zeta\in[k-1+\tau,k+\tau]$ or $\zeta\in[0,k+\tau]$, respectively). By relative compactness of $\mathcal I$ we know that $\|\mathcal X(s,\tilde\lambda)\|\lesssim 1$ uniformly
in $s$ and $\tilde\lambda\in\mathcal I$. Thus, integration over $\tau$ gives the desired bounds \eqref{eq:4.17},
\begin{equation}
   \| \mathcal Y_j(k+1,k,\xi) - \mathcal X( \tilde\lambda)\|\le\int_0^1 \|\mathcal X(\tau,\tilde\lambda)\| \|A_j(k+\tau,\xi)-B(\tau,\tilde\lambda)\|\d\tau\lesssim \delta_j\frac{\lambda(t_j)}{\Lambda(t_j)} 
\end{equation}
uniform in $k$, $j$ and $\xi$ and using $\|\mathcal Y_j(k+\tau,k,\xi) \|\lesssim 1$, $\tau\in[0,1]$,
as consequence of $\nu_j \frac{\lambda(t_j)}{\Lambda(t_j)}\lesssim1$ also \eqref{eq:4.16},
\begin{align}
   &\| \mathcal Y_j(k+1,k,\xi) -\mathcal Y_j(k,k-1,\xi)\|\notag\\&\le\int_0^1 \|\mathcal Y_j(k-1+\tau,k-1,\xi)\| \|A_j(k+\tau,\xi)-A_j(k-1+\tau,\xi)\|\d\tau\lesssim\frac{\delta_j}{\nu_j}\frac{\lambda(t_j)}{\Lambda(t_j)}
\end{align}
uniform in $k$, $j$ and $\xi$.

{\sl Step 2.} In a second step we want to compare $\mathcal Y_j(\nu_j,0,\xi)$ with $\mathcal X^{\nu_j}(\tilde\lambda)$.  For this we denote by $M_{j,k}(\xi)$ diagonaliser of $\mathcal Y_j(k,k-1,\xi)$ 
and $M(\tilde\lambda)$ of $\mathcal X(\tilde\lambda)$ which are of bounded condition uniform in $j$ and close to each other. Furthermore, we denote by $D_{j,k}(\xi)$ and $D(\tilde\lambda)$ the corresponding diagonal matrices (having the big eigenvalue as upper left corner entry). Then
\begin{align}\label{eq:4.22}
&M^{-1}(\tilde\lambda)\mathcal Y_j(\nu_j,0,\xi)M(\tilde\lambda)\notag \\ 
&=M^{-1}(\tilde\lambda)\mathcal Y_j(\nu_j,\nu_j-1,\xi)\cdots\mathcal Y_j(2,1,\xi) \mathcal Y_j(1,0,\xi)M(\tilde\lambda)\notag\\
& = M^{-1}(\tilde\lambda)M_{j,\nu_j}(\xi) D_{j,\nu_j}(\xi) M_{j,\nu_j}^{-1}(\xi)M_{j,\nu_j-1}(\xi) \cdots 
D_{j,1}(\xi) M_{j,1}^{-1}(\xi)M(\tilde\lambda)  \notag\\
&=  (\mathrm I+G_{j,\nu_j+1}(\xi))D(\tilde\lambda)(\mathrm I+G_{j,\nu_j}(\xi))  \cdots D(\tilde\lambda) (\mathrm I+G_{j,2}(\xi)) D(\tilde\lambda)  (\mathrm I+G_{j,1}(\xi)) ,
\end{align}
where $G_{j,k}(\xi)  = D^{-1}(\tilde\lambda) D_{j,k}(\xi)M_{j,k}^{-1}(\xi)M_{j,k-1}(\xi)-\mathrm I$ and for convenience $M_{j,0}(\xi)=M(\tilde\lambda)=M_{j,\nu_j+1}(\xi)$, $D_{j,\nu_j+1}(\xi)=D(\tilde\lambda)$.

We need to look at the diagonaliser in more detail. Due to Liouville theorem we know that $\det\mathcal Y_j(k,k-1,\xi)=\det\mathcal X(\tilde\lambda)=1$ and the matrices $M_{j,k}(\xi)$ may be expressed in terms of the entries of $\mathcal Y_j(k,k-1,\xi)$ and their eigenvalues. If we denote them as $y^{(j,k)}_{mn}(\xi)$ and the eigenvalues as $\mu^{\pm1}_{j,k}(\xi)$ and assume for simplicity that $|y^{(j,k)}_{11}| \le |y^{(j,k)}_{22}|$ 
a suitable diagonaliser is
\begin{equation}
   M_{j,k}(\xi) = \begin{pmatrix} \frac{y^{(j,k)}_{21}}{\mu_{j,k}^{-1}-y^{(j,k)}_{22}}&1\\1& \frac{y^{(j,k)}_{12}}{\mu_{j,k}-y^{(j,k)}_{11}}\end{pmatrix}.
\end{equation}
A similar formula holds for $M(\tilde\lambda)$. Due to the estimates of Step 1 a short calculation implies
that $M_{j,k}^{-1}(\xi)M_{j,k-1}(\xi)$ approximates the identity,
\begin{equation}
   \|M_{j,k}^{-1}(\xi)M_{j,k-1}(\xi)-\mathrm I\| \lesssim \frac{\delta_j}{\nu_j}\frac{\lambda(t_j)}{\Lambda(t_j)}
\end{equation}
unform in $j$, and therefore
\begin{align}
  &\|G_{j,k}(\xi)\| \lesssim \frac{\delta_j}{\nu_j} \frac{\lambda(t_j)}{\Lambda(t_j)},\qquad k=2,\dots,\nu_j,\\
  \intertext{and similarly}
  &\|G_{j,1}(\xi)\|, \|G_{j,\nu_j+1}(\xi)\|\lesssim \delta_j \frac{\lambda(t_j)}{\Lambda(t_j)}\lesssim 1.
\end{align}
Therefore, the right hand side of \eqref{eq:4.22} can be written as $D^{\nu_j}(\tilde\lambda)$ plus a remainder of size 
\begin{equation}
\lesssim   \sum_{k=1}^{\nu_j-1} \binom{\nu_j-1}{k} \mu^{\nu_j-k} \left( \frac{\delta_j}{\nu_j}\frac{\lambda(t_j)}{\Lambda(t_j)}\right)^k 
   = \mu^{\nu_j} \left(\left(1+\frac1\mu  \frac{\delta_j}{\nu_j}\frac{\lambda(t_j)}{\Lambda(t_j)}\right)^{\nu_j-1} -1 \right)
\end{equation}
uniformly in $j$ and with $\| D(\tilde\lambda)\|\sim \mu$. Due to our assumption $\delta_j\lambda(t_j)/\Lambda(t_j)\to0$ and therefore the expression in brackets behaves like $\exp(\mu^{-1} \frac{\nu_j-1}{\nu_j}\delta_j\lambda(t_j)/\Lambda(t_j))-1$, which tends to zero as $j$ approaches infinity. 
Choosing $j$ large enough to bound this expression by $1/2$ and application of Bauer-Fike theorem proves the desired statement. 
$\Box$\end{pf}

\subsubsection{Choice of sequences} We assume $j$ is large enough. Then the eigenvalues of $\mathcal Y_j(\nu_j,0,\xi)$ are distinct and we are allowed to choose $V_j(\xi)\in C_0^\infty (\Omega_j)$  in such a way that it is normalised in $L^2$-sense and coincides for each fixed $\xi$ on its support with an eigenvector of  $\mathcal Y_j(\nu_j,0,\xi)$ corresponding to the large eigenvalue. Then we solve $\D_t V = A(t,\xi)V$ with $V(t_j,\xi)=V_j(\xi)$ and denote by $u_j$ the corresponding solution of the original problem. This yields a sequence of solutions with a remarkable property.
As consequence of Lemma~\ref{lem:4.2} we obtain uniformly in $j$, $j$ large,
\begin{equation}\label{eq:4.28}
   \mathbb E_\lambda(u_j;t_j+\delta_j) \gtrsim \mu^{2\nu_j}\mathbb E_\lambda(u_j;t_j)= \mu^{2\nu_j}.
\end{equation}
This estimate contradicts with the estimate of Lemma~\ref{lem:2.5}, which implies uniform in 
$j$, $j$ large,
\begin{equation}\label{eq:4.29}
   \mathbb E_\lambda(u_j;t_j+\delta_j) \lesssim \frac{\lambda(t_j)}{\lambda(t_j+\delta_j)} \mathbb E_\lambda(u_j;t_j)\approx 1,
\end{equation}
provided that 
$\Theta(t_j)|\xi|\approx \frac{\nu_j}{\delta_j} \frac{\Theta(t_j)}{\lambda(t_j)}\to\infty$, i.e., $[t_j,t_j+\delta_j]\times\Omega_j$ belongs to the hyperbolic zone for large $j$. The estimates \eqref{eq:4.28} and \eqref{eq:4.29} contradict each other. 

Thus, if we manage to construct sequences $t_j$, $\delta_j$ and $\nu_j$ such that all requirements are satisfied, a counter-example is found. We will do this for all our example classes. 

\begin{expl}\label{expl1-3}
(Counter-example, polynomial case) 
Let $\lambda(t)=(1+t)^p$ for some $p\ge0$ and $\Theta(t)=(1+t)^q$, $-1\le q<p$. We construct admissible sequences such that (A1)--(A3) hold, but (A4'') is violated in the sense that such an estimate holds only for a given arbitrarily small {\em negative} exponent. 

We choose $t_j=2^j$, $\delta_j=2^{j(q-p+1)-1}$  and $\nu_j = \lceil 2^{j\epsilon(p-q)} \rceil$. Stabilisation is ensured and (A1) -- (A3) are valid.  By construction $\lambda(t_j+\delta_j) \approx \lambda(t_j)$
and $\Lambda(t_j+\delta_j)\approx\Lambda(t_j)$ holds uniformly in $j$ and $\delta_j \lambda(t_j)/\Lambda(t_j)\approx 2^{j (q-p)}\to0$. Thus, Lemma~\ref{lem:4.2} can be applied. It remains to check the geometry restriction arising from the zone. It follows on $\Omega_j$ that $\Theta(t_j)|\xi|\approx \frac{\nu_j\Theta(t_j)}{\delta_j\lambda(t_j)}\approx 2^{j\epsilon(p-q)}\to\infty$. 

We check how closely (A4'') is violated. Since derivatives behave like multiplications with $\nu_j/\delta_j$ on the interval $[t_j,t_j+\delta_j]$, the best possible choice of $\Xi(t)$ would be 
\begin{equation}\label{eq:4.30}
\Xi(t) = (1+t)^{q-p+1-\epsilon(p-q)}=\left(\frac{\lambda(t)}{\Theta(t)}\left(\frac{\Theta(t)}{\Lambda(t)}\right)^{-\epsilon}\right)^{-1}
\end{equation}
in contrast to (A4'').
\end{expl}

\begin{expl}\label{expl2-3}
(Counter-example, supra-polynomial case)
Let $\lambda(t)=\exp(t^\alpha)$ with $\alpha\in(0,1)$ and $\Theta(t)=t^{-\beta}\exp(t^\alpha)$, $\beta\ge\alpha-1$. We choose $t_j=j^{1/\alpha}$ and $\delta_j=j^{-\beta/\alpha}$, such (A1)--(A3) are valid.
Furthermore, we choose $\nu_j$ as $\nu_j = \lceil j^{\epsilon(\beta-\alpha+1)/\alpha}\rceil$, $\epsilon>0$. 
It is evident that $\lambda(t_j+\delta_j)\approx \lambda(t_j)$ holds and similarly for the primitive. Again by construction $\delta_j\lambda(t_j)/\Lambda(t_j)\approx j^{-(\beta-\alpha+1)/\alpha}$
tends to zero if $\beta>\alpha-1$ such that Lemma~\ref{lem:4.2} applies. Furthermore, 
$ \frac{\nu_j\Theta(t_j)}{\delta_j\lambda(t_j)}\approx j^{-\epsilon(\beta-\alpha+1)/\alpha}\to\infty$ and the counter-example is constructed. 
Derivatives behave like multiplication with $(1+t)^{\beta+\epsilon(\beta-\alpha+1)}$, i.e. (A4'') with exponent $-\epsilon$ (cf. equation \eqref{eq:4.30}). 
\end{expl}

\begin{expl}\label{expl3-3}
(Counter-example, exponential case)
Let $\lambda(t)=\mathrm e^t$ and $\Theta(t)=\mathrm e^{at}$, $a<1$. We choosing $t_j=j$, 
$\delta_j=\mathrm e^{j(a-1)}$ and $\nu_j=\lceil \mathrm e^{j\epsilon(1-a)}\rceil$, $\epsilon>0$. Then (A1)--(A3) hold.
From $\delta_j\to0$ we conclude $\lambda(t_j+\delta_j)\approx\lambda(t_j)$, the primitive is the same function. Furthermore, $\delta_j\lambda(t_j)/\Lambda(t_j)\approx \delta_j\to0$ and Lemma~\ref{lem:4.2} applies and the geometry restriction $\frac{\nu_j\Theta(t_j)}{\delta_j\lambda(t_j)}\approx \mathrm e^{j\epsilon(1-a)}\to\infty$ is valid.
The behaviour of derivatives is described by $\Xi(t)=\mathrm e^{(a-1-\epsilon(1-a))t}$, thus again (A4'') holds only with exponent $-\epsilon$ (cf. equation \eqref{eq:4.30}). 
\end{expl}

Hence, in all cases there exists a coefficient function $a(t)$ satisfying (A1) -- (A3) and violating (A4'') to arbitrary small order for which the statement of Lemma~\ref{lem:2.5} is false. 

We finally want to discuss how to conclude a counter-example for the estimate of Theorem~\ref{thm:3.2}. For this we use the same idea as above, but estimate the corresponding Cauchy data on the level $t=0$.  
Let for this $V_j$ and $u_j$ be constructed as above, $\mathbb E(u_j;t_j)=1$ and $\mathbb E(u_j;t_j+\delta_j)\gtrsim \mu^{2\nu_j}$ uniform in the sequence $u_j$ and assume that (A1)--(A3) hold true. 
We are going to estimate 
\begin{equation} 
   \mathcal E(0,t_j,\xi) V_j(\xi) = \mathcal E(0,t_\xi^{(1)},\xi) \mathcal E(t_\xi^{(1)},t_\xi^{(2)},\xi) \mathcal E(t_\xi^{(2)},t_j,\xi)V_j(\xi).
\end{equation}
The first two factors satisfy Lemma~\ref{lem:2.2} and~\ref{lem:2.8}, respectively. Both Lemmata are true as consequence of the above assumptions. For the third one we use 
\begin{equation}
   \|\mathcal E(t_\ell+\delta_\ell,t_{\ell+1},\xi)\| \approx \frac{\sqrt{\lambda(t_{\ell+1})}}{\sqrt{\lambda(t_\ell+\delta_\ell)}}
\end{equation}
uniform in $\ell$ with $\Lambda(t_\ell)|\xi|\ge N$ as consequence of Lemma~\ref{lem:2.1} in combination with 
\begin{equation}
   \|\mathcal E(t_\ell,t_\ell+\delta_\ell,\xi)\| \lesssim \mathrm e^{c\nu_\ell}, \qquad  c=\sup_\tau \frac{|b'(\tau)|}{1+b(\tau)},
\end{equation}
following from Gronwall inequality. Combining all these estimates we get for the Cauchy data 
$u_{j,1}$ and $u_{j,2}$ corresponding to the solution $u_j$
\begin{equation}
  \|u_{j,1}\|_{H^1} + \|u_{j,2}\|_{L^2} \lesssim S^{j/2}\frac{t_{\xi}^{(1)} \sqrt{\lambda(t_\xi^{(1)})}}{\sqrt{\lambda(t_j)}} \exp\left(c \sum_{\ell=\ell_0}^{j-1} \nu_\ell \right),
\end{equation}
where $S=\sup_j{\lambda(t_j+\delta_j)}/{\lambda(t_j)}$.
Using (A1$_+$) and the definition of $\mathcal I_j$ it follows that $t_{\xi}^{(1)} \sqrt{\lambda(t_\xi^{(1)})}\lesssim \Lambda(t_\xi^{(1)})\approx|\xi|^{-1}\approx \delta_j\lambda(t_j)/\nu_j$.  If Theorem~\ref{thm:3.2} would be true, it would imply
\begin{equation}
\mu^{2\nu_j} \lesssim  \mathbb E_\lambda(t_j+\delta_j;u_j) \lesssim \frac{\delta_j^2}{\nu_j^2} \lambda^2(t_j) S^j \exp\left(2c \sum_{\ell=\ell_0}^{j-1} \nu_\ell \right).
\end{equation}
This gives a contradiction if 
\begin{equation}\label{eq:4.35}
 \frac{\delta_j^2}{\nu_j^2} \lambda^2(t_j) S^j \exp\left(2c \sum_{\ell=\ell_0}^{j-1} \nu_\ell  - 2\nu_j \log\mu \right)\to 0, \qquad j\to\infty.
\end{equation}

We are going to check this for the previously constructed counter-example in the polynomial case. 

\begin{expl}(Counter-example, polynomial case)
We follow Example~\ref{expl1-3} for $\lambda(t)=(1+t)^p$, $\Theta(t)=(1+t)^q$, however with a minor change. We choose sequences $t_j=\sigma^j$, $\delta_j=\sigma^{j(q-p+1)-1}$ and $\nu_j=\lceil\sigma^{j\epsilon(p-q)}\rceil$, $\epsilon>0$. All the previous considerations and conditions transfer, thus choosing $\epsilon$ small enough will closely violate (A4''). If we now consider the condition \eqref{eq:4.35},
the first factors increase exponentially like $S^j\sigma^{j(q+1-\epsilon(p-q))}$, while the second exponential can be estimated by
\begin{equation}
   \exp\left(2c \frac{\sigma^{j\epsilon(p-q)}-1}{\sigma^{\epsilon(p-q)}-1} - 2 \sigma^{j\epsilon(p-q)} \log\mu \right) \lesssim \exp(-c' \sigma^{j\epsilon(p-q)}), \qquad c'>0,
\end{equation}
provided $\sigma$ is chosen large, ${c}/{(\sigma^{\epsilon(p-q)}-1)}<\log\mu$. Thus we obtain a counter-example to Theorem~\ref{thm:3.2}.
\end{expl}

\begin{expl}(Counter-example, exponential case) We follow Example~\ref{expl3-3} with $\lambda(t)=\mathrm e^t$, $\Theta(t)=\mathrm e^{at}$  and choose the sequences $t_j=\sigma j$, $\delta_j =\mathrm e^{\sigma j(a-1)}$ and $\nu_j=\lceil\mathrm e^{\sigma j\epsilon(1-a)}\rceil$ with a new additional parameter $\sigma$. For any choice of $\sigma>0$ the reasoning of Example~\ref{expl3-3} remains true. Furthermore,  \eqref{eq:4.35} follows provided that $\sigma$ is chosen big enough, i.e. if $c / (e^{\sigma \epsilon(1-a)}-1) < \log\mu$ holds. 
\end{expl}


\end{document}